\documentclass[reqno]{amsart}
\usepackage{amssymb,hyperref}
\usepackage[dvips]{graphicx}

\makeatletter

\@addtoreset{equation}{section}

\theoremstyle{plain}
\newtheorem{thm}{Theorem}[section]
\newtheorem{lem}{Lemma}[section]
\newtheorem{cor}{Corollary}[section]
\newtheorem{ex}{Example}[section]

\theoremstyle{remark}
\newtheorem{remark}{Remark}[section]

\title[Coefficient conditions for harmonic close-to-convex functions]
{Coefficient conditions for\\
harmonic close-to-convex functions}

\author{Toshio Hayami}
\address{Toshio Hayami \newline
Department of Mathematics, \newline
Kinki University \newline
Higashi-Osaka, Osaka 577-8502, \newline
Japan}
\email{ha\_ya\_to112@hotmail.com}

\subjclass[2010]{Primary 30C45, Secondary 58E20.}
\keywords{Coefficient condition, harmonic function, univalent function, close-to-convex function.}
\date{}

\begin{document}

\begin{abstract}
New sufficient conditions, concerned with the coefficients of harmonic functions $f(z)=h(z)+\overline{g(z)}$ in the open unit disk $\mathbb{U}$ normalized by $f(0)=h(0)=h'(0)-1=0$, for $f(z)$ to be harmonic close-to-convex functions are discussed. Furthermore, several illustrative examples and the image domains of harmonic close-to-convex functions satisfying the obtained conditions are enumerated.
\end{abstract}

\begin{flushleft}
This paper was published in the journal: \\
Abstr. Appl. Anal. Vol.2012, Article ID 413965, 12 pages.\\
\url{http://www.hindawi.com/journals/aaa/2012/413965/}
\end{flushleft}
\hrule

\

\

\maketitle

\section{Introduction}

\

For a continuous complex-valued function $f(z)=u(x,y)+iv(x,y)$\ \ $(z=x+iy)$, we say that $f(z)$ is harmonic in the open unit disk $\mathbb{U}=\{z\in \mathbb{C}:|z|<1\}$ if both $u(x,y)$ and $v(x,y)$ are real harmonic in $\mathbb{U}$, that is, $u(x,y)$ and $v(x,y)$ satisfy the Laplace's equations
$$
\Delta u=u_{xx}+u_{yy}=0\qquad {\rm and}\qquad \Delta v=v_{xx}+v_{yy}=0.
$$
A complex-valued harmonic function $f(z)$ in $\mathbb{U}$ is given by $f(z)=h(z)+\overline{g(z)}$ where $h(z)$ and $g(z)$ are analytic in $\mathbb{U}$. We call $h(z)$ and $g(z)$ the analytic part and the co-analytic part of $f(z)$, respectively. A necessary and sufficient condition for $f(z)$ to be locally univalent and sense-preserving in $\mathbb{U}$ is $|h'(z)|>|g'(z)|$ in $\mathbb{U}$ (see, \cite{CS} or \cite{L}). Let $\mathcal{H}$ denote the class of harmonic functions $f(z)$ in $\mathbb{U}$ with $f(0)=h(0)=0$ and $h'(0)=1$. Thus, every normalized harmonic function $f(z)$ can be written by
$$
f(z)=h(z)+\overline{g(z)}=z+\sum\limits_{n=2}^{\infty}a_n z^n+\overline{\sum\limits_{n=1}^{\infty}b_n z^n}\in \mathcal{H}
$$
where $a_1=1$ and $b_0=0$, for convenience.

\

We next denote by $\mathcal{S}_{\mathcal{H}}$ the class of functions $f(z)\in \mathcal{H}$ which are univalent and sense-preserving in $\mathbb{U}$. Since the sense-preserving property of $f(z)$, we see that $|b_1|=|g'(0)|<|h'(0)|=1$. If $g(z)\equiv 0$, then $\mathcal{S}_{\mathcal{H}}$ reduces to the class $\mathcal{S}$ consisting of normalized analytic univalent functions. Furthermore, for every function $f(z)\in \mathcal{S}_{\mathcal{H}}$, the function

\

$$
F(z)=\dfrac{f(z)-\overline{b_1 f(z)}}{1-|b_1|^2}=z+\sum\limits_{n=2}^{\infty}\dfrac{a_n-\overline{b_1}b_n}{1-|b_1|^2}z^n+\overline{\sum\limits_{n=2}^{\infty}\dfrac{b_n-b_1 a_n}{1-|b_1|^2}z^n}
$$
is also a member of $\mathcal{S}_{\mathcal{H}}$. Therefore, we consider the subclass $\mathcal{S}_{\mathcal{H}}^{0}$ of $\mathcal{S}_{\mathcal{H}}$ defined as
$$
\mathcal{S}_{\mathcal{H}}^{0}=\left\{f(z)\in \mathcal{S}_{\mathcal{H}}:b_1=g'(0)=0\right\}.
$$
Conversely, if $F(z)\in \mathcal{S}_{\mathcal{H}}^{0}$, then $f(z)=F(z)+\overline{b_1 F(z)}\in \mathcal{S}_{\mathcal{H}}$ for any $b_1$ $(|b_1|<1)$.

\

We say that a domain $\mathbb{D}$ is a close-to-convex domain if the complement of $\mathbb{D}$ can be written as a union of non-intersecting half-lines (except that the origin of one half-line may lie on one of the other half-lines). Let $\mathcal{C}$, $\mathcal{C}_{\mathcal{H}}$ and $\mathcal{C}_{\mathcal{H}}^{0}$ be the respective subclasses of $\mathcal{S}$, $\mathcal{S}_{\mathcal{H}}$ and $\mathcal{S}_{\mathcal{H}}^{0}$ consisting of all functions $f(z)$ which map $\mathbb{U}$ onto a certain close-to-convex domain.\\

Bshouty and Lyzzaik \cite{BL} have stated the following result.

\begin{thm} \label{thm1}
If $f(z)=h(z)+\overline{g(z)}\in \mathcal{H}$ satisfies
$$
g'(z)=zh'(z)\qquad and\qquad {\rm Re}\left(1+\dfrac{zh''(z)}{h'(z)}\right)>-\dfrac{1}{2}
$$
for all $z\in \mathbb{U}$, then $f(z)\in \mathcal{C}_{\mathcal{H}}^{0}\subset \mathcal{S}_{\mathcal{H}}^{0}$.
\end{thm}

\

A simple and interesting example is below.

\begin{ex} \label{ex1}
The function
$$
f(z)=\dfrac{1-(1-z)^2}{2(1-z)^2}+\overline{\dfrac{z^2}{2(1-z)^2}}=z+\sum\limits_{n=2}^{\infty}\dfrac{n+1}{2}z^n+\sum\limits_{n=2}^{\infty}\dfrac{n-1}{2}\overline{z}^n
$$
satisfies the conditions of Theorem \ref{thm1}, and therefore $f(z)$ belongs to the class $\mathcal{C}_{\mathcal{H}}^{0}$. We now show that $f(\mathbb{U})$ is actually close-to-convex domain. It follows that
\begin{eqnarray*}
f(z) &=& \left(\dfrac{z}{2(1-z)^2}+\dfrac{z}{2(1-z)}\right)+\overline{\left(\dfrac{z}{2(1-z)^2}-\dfrac{z}{2(1-z)}\right)}\\
 & & \\
 &=& {\rm Re}\left(\dfrac{z}{(1-z)^2}\right)+i{\rm Im}\left(\dfrac{z}{1-z}\right).
\end{eqnarray*}
Setting
$$
f(re^{i\theta})=\dfrac{-2r^2+r(1+r^2)\cos \theta}{(1+r^2-2r\cos \theta)^2}+\dfrac{r\sin \theta}{1+r^2-2r\cos \theta}i=u+iv
$$
for any $z=re^{i\theta}\in \mathbb{U}$ $(0\leqq r<1,\ 0\leqq \theta<2\pi)$, we see that
$$
-4(u+v^2)=\dfrac{4r(r-\cos \theta)(1-r\cos \theta)}{(1+r^2-2r\cos \theta)^2}=\dfrac{4r(r-t)(1-rt)}{(1+r^2-2rt)^2}\equiv \phi(t)\qquad (-1\leqq t=\cos \theta\leqq 1).
$$
Since
$$
\phi'(t)=\dfrac{-4r(1-r^2)^2}{(1+r^2-2rt)^3}\leqq 0,
$$
we obtain that
$$
\phi(t)\leqq \phi(-1)=\dfrac{4r}{(1+r)^2}\equiv \psi(r).
$$
Also, noting that
$$
\psi'(r)=\dfrac{4(1-r)}{(1+r)^3}>0,
$$
we know that
$$
\psi(r)<\psi(1)=1
$$
which implies that
$$
u>-v^2-\dfrac{1}{4}.
$$
Thus, $f(z)$ maps $\mathbb{U}$ onto the following close-to-convex domain.

\begin{figure*}[h]
\begin{center}
\includegraphics[width=5cm,clip]{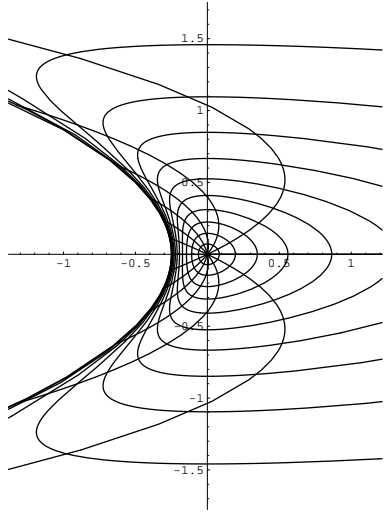}
\end{center}
\caption{The image of $f(z)=\dfrac{1-(1-z)^2}{2(1-z)^2}+\overline{\dfrac{z^2}{2(1-z)^2}}$.}
\end{figure*}
\end{ex}

\

\begin{remark}
Let $\mathcal{M}$ be the class of all functions satisfying the conditions of Theorem \ref{thm1}. Then, it was earlier conjectured by Mocanu \cite{M1,M2} that $\mathcal{M}\subset \mathcal{S}_{\mathcal{H}}^{0}$. Furthermore, we can immediately see that the function $f(z)$ in Example \ref{ex1} is a member of the class $\mathcal{M}$ and it shows that $f(z)\in \mathcal{M}$ is not necessarily starlike with respect to the origin in $\mathbb{U}$ ($f(z)$ is starlike with respect to the origin in $\mathbb{U}$ if and only if $tw\in f(\mathbb{U})$ for all $w\in f(\mathbb{U})$ and $t$ $(0\leqq t\leqq 1)$).
\end{remark}

\

\begin{remark}
For the function $f(z)=h(z)+\overline{g(z)}\in \mathcal{H}$ given by
$$
g'(z)=z^{n-1}h'(z)\qquad (n=2,3,4,\cdots),
$$
letting $w(t)=f(e^{it})=h(e^{it})+\overline{g(e^{it})}$\ \ $(-\pi\leqq t<\pi)$, we know that
$$
{\rm Im}\left(\dfrac{w''(t)}{w'(t)}\right)\leqq 0\qquad (-\pi\leqq t<\pi)
$$
which means that $f(z)$ maps the unit circle $\partial \mathbb{U}=\left\{z\in \mathbb{C}:|z|=1\right\}$ onto a union of several concave curves (see, \cite[Theorem 2.1]{HO}).
\end{remark}

\

Jahangiri and Silverman \cite{JS} have given the following coefficient inequality for $f(z)\in \mathcal{H}$ to be in the class $\mathcal{C}_{\mathcal{H}}$.

\begin{thm} \label{thm2}
If $f(z)\in \mathcal{H}$ satisfies
$$
\sum\limits_{n=2}^{\infty}n|a_n|+\sum\limits_{n=1}^{\infty}n|b_n|\leqq 1,
$$
then $f(z)\in \mathcal{C}_{\mathcal{H}}$.
\end{thm}

\

\begin{ex}
The function
$$
f(z)=z+\dfrac{1}{5}\overline{z}^5
$$
belongs to the class $\mathcal{C}_{\mathcal{H}}^{0}\subset \mathcal{C}_{\mathcal{H}}$ and satisfies the condition of Theorem \ref{thm2}. Indeed, $f(z)$ maps $\mathbb{U}$ onto the following hypocycloid of six cusps (cf. \cite{D} or \cite{HO}).

\

\begin{figure*}[h]
\begin{center}
\includegraphics[width=7cm,clip]{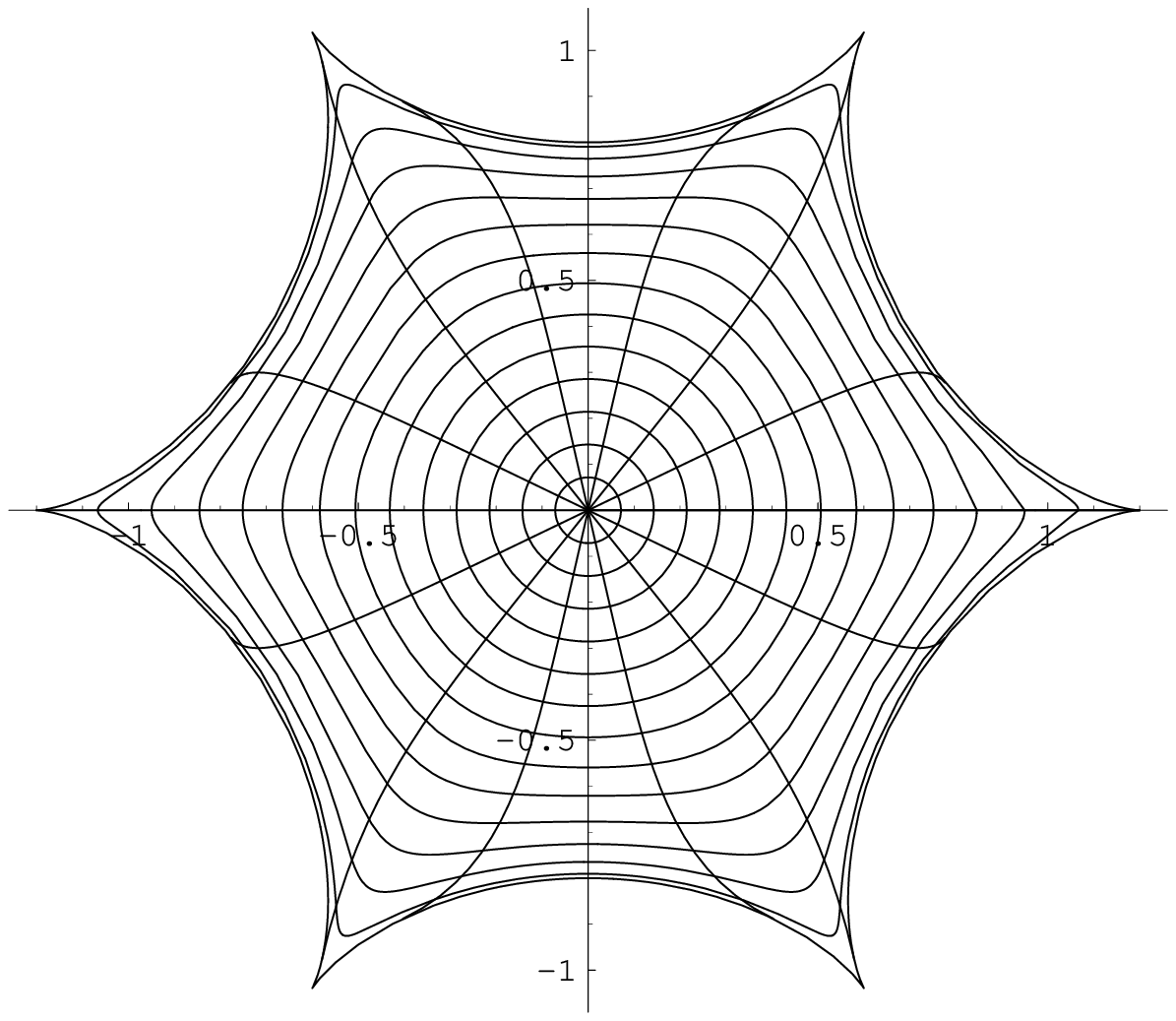}
\end{center}
\caption{The image of $f(z)=z+\dfrac{1}{5}\overline{z}^5$.}
\end{figure*}
\end{ex}

\

The object of this paper is to find some sufficient conditions for functions $f(z)\in \mathcal{H}$ to be in the class $\mathcal{C}_{\mathcal{H}}$. In order to establish our results, we have to recall here the following lemmas due to Clunie and Sheil-small \cite{CS}.

\begin{lem} \label{lem1}
If $h(z)$ and $g(z)$ are analytic in $\mathbb{U}$ with $|h'(0)|>|g'(0)|$ and $h(z)+\varepsilon g(z)$ is close-to-convex for each $\varepsilon$ $(|\varepsilon|=1)$, then $f(z)=h(z)+\overline{g(z)}$ is harmonic close-to-convex.
\end{lem}

\

\begin{lem} \label{lem2}
If $f(z)=h(z)+\overline{g(z)}$ is locally univalent in $\mathbb{U}$ and $h(z)+\varepsilon g(z)$ is convex for some $\varepsilon$ $(|\varepsilon|\leqq 1)$, then $f(z)$ is univalent close-to-convex.
\end{lem}

\

We also need the following result due to Hayami, Owa and Srivastava \cite{HOS}.

\begin{lem} \label{lem3}
If a function $F(z)=z+\sum\limits_{n=2}^{\infty}A_n z^n$ is analytic in $\mathbb{U}$ and satisfies\\
\\
${\displaystyle 
\sum\limits_{n=2}^{\infty}\left[\left|\sum\limits_{k=1}^{n}\left\{\sum\limits_{j=1}^{k}(-1)^{k-j}j(j+1)\left(
\begin{array}{c}
\alpha\\
k-j
\end{array}
\right)A_j\right\}\left(
\begin{array}{c}
\beta\\
n-k
\end{array}
\right)\right|\right.
}$
\begin{flushright}
${\displaystyle 
+\left.\left|\sum\limits_{k=1}^{n}\left\{\sum\limits_{j=1}^{k}(-1)^{k-j}j(j-1)\left(
\begin{array}{c}
\alpha\\
k-j
\end{array}
\right)A_j\right\}\left(
\begin{array}{c}
\beta\\
n-k
\end{array}
\right)\right|\right]\leqq 2
}$
\end{flushright}
for some real numbers $\alpha$ and $\beta$, then $F(z)$ is convex in $\mathbb{U}$.
\end{lem}

\

\section{Main results}

\

Our first result is contained in

\begin{thm} \label{thm3}
If $f(z)\in \mathcal{H}$ satisfies the following condition
$$
\sum\limits_{n=2}^{\infty}\left|na_n-e^{i\varphi}(n-1)a_{n-1}\right|+\sum\limits_{n=1}^{\infty}\left|nb_n-e^{i\varphi}(n-1)b_{n-1}\right|\leqq 1
$$
for some real number $\varphi$ $(0\leqq \varphi<2\pi)$, then $f(z)\in \mathcal{C}_{\mathcal{H}}$.
\end{thm}

\

\begin{proof}
Let $F(z)=z+\sum\limits_{n=2}^{\infty}A_n z^n$ be analytic in $\mathbb{U}$. If $F(z)$ satisfies
$$
\sum\limits_{n=2}^{\infty}\left|nA_n-e^{i\varphi}(n-1)A_{n-1}\right|\leqq 1
$$
then it follows that
\begin{eqnarray*}
\left|(1-e^{i\varphi}z)F'(z)-1\right| &=& \left|\sum\limits_{n=2}^{\infty}\left(nA_n-e^{i\varphi}(n-1)A_{n-1}\right)z^{n-1}\right|\\
 & & \\
 &\leqq& \sum\limits_{n=2}^{\infty}\left|nA_n-e^{i\varphi}(n-1)A_{n-1}\right|\cdot |z|^{n-1}\\
 & & \\
 &<& \sum\limits_{n=2}^{\infty}\left|nA_n-e^{i\varphi}(n-1)A_{n-1}\right|\ \leqq\ 1\qquad (z\in \mathbb{U}).
\end{eqnarray*}
This gives us that
$$
{\rm Re}\left((1-e^{i\varphi}z)F'(z)\right)>0\qquad (z\in \mathbb{U}),
$$
that is, that $F(z)\in \mathcal{C}$. Then, it is sufficient to prove that
$$
F(z)=\dfrac{h(z)+\varepsilon g(z)}{1+\varepsilon b_1}=z+\sum\limits_{n=2}^{\infty}\dfrac{a_n+\varepsilon b_n}{1+\varepsilon b_1}z^n\in \mathcal{C}
$$
for each $\varepsilon$ $(|\varepsilon|=1)$ by Lemma \ref{lem1}. From the assumption of the theorem, we obtain that\\
\\
${\displaystyle 
\sum\limits_{n=2}^{\infty}\left|n\dfrac{a_n+\varepsilon b_n}{1+\varepsilon b_1}-e^{i\varphi}(n-1)\dfrac{a_{n-1}+\varepsilon b_{n-1}}{1+\varepsilon b_1}\right|
}$
$$
\leqq \dfrac{1}{1-|b_1|}\sum\limits_{n=2}^{\infty}\left[\left|na_n-e^{i\varphi}(n-1)a_{n-1}\right|+\left|nb_n-e^{i\varphi}(n-1)b_{n-1}\right|\right]\leqq \dfrac{1-|b_1|}{1-|b_1|}=1.
$$
This completes the proof of the theorem.
\end{proof}

\

\begin{ex}
The function
$$
f(z)=-\log(1-z)+\overline{\Bigl(-mz-\log(1-z)\Bigr)}=z+\sum\limits_{n=2}^{\infty}\dfrac{1}{n}z^n+(1-m)\overline{z}+\sum\limits_{n=2}^{\infty}\dfrac{1}{n}\overline{z}^n\quad (0<m\leqq 1)
$$
satisfies the condition of Theorem \ref{thm3} with $\varphi=0$ and belongs to the class $\mathcal{C}_{\mathcal{H}}$. In particular, putting $m=1$, we obtain the following.

\

\begin{figure*}[h]
\begin{center}
\includegraphics[width=5cm,clip]{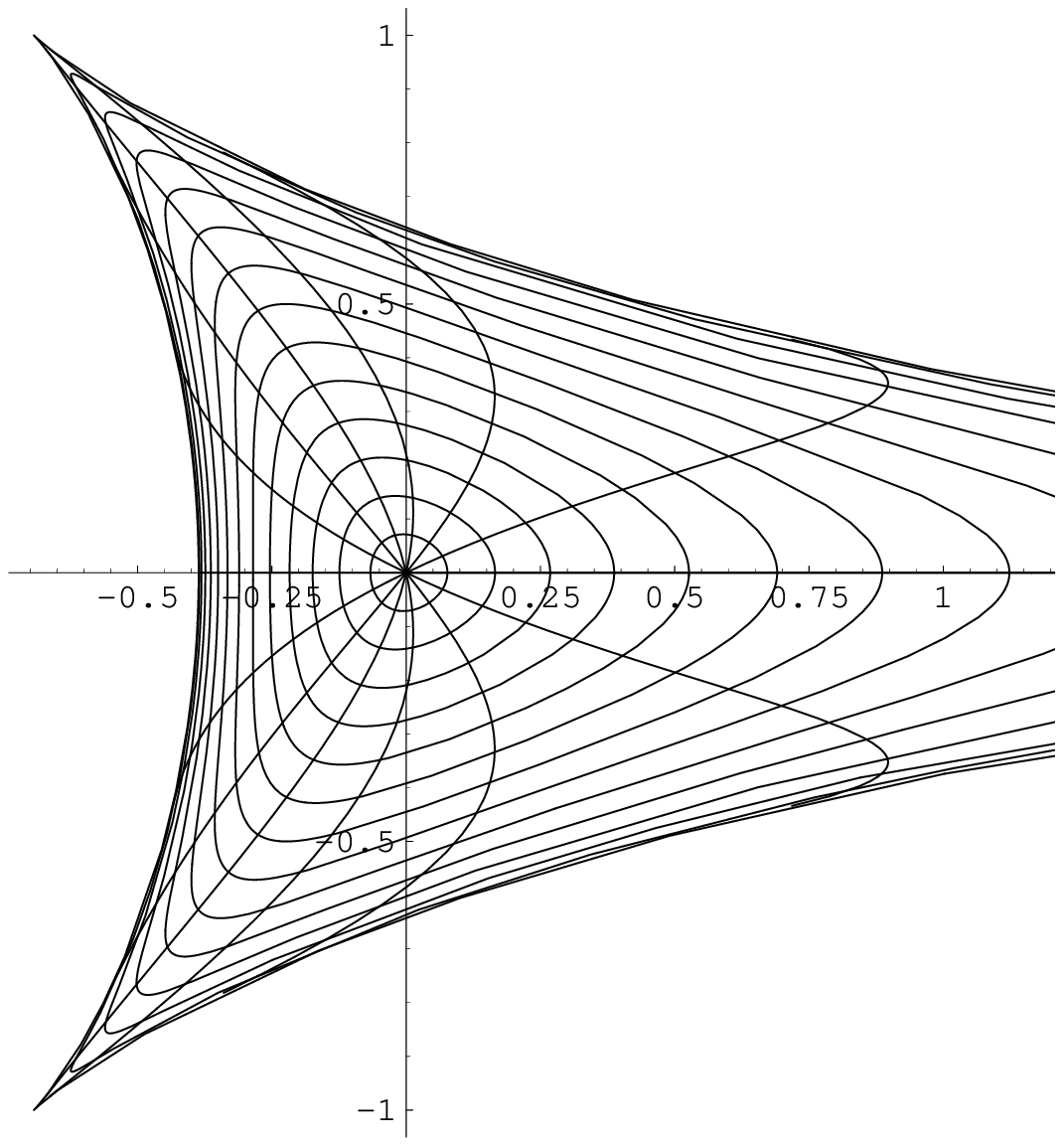}
\end{center}
\caption{The image of $f(z)=-\overline{z}-2\log|1-z|$.}
\end{figure*}
\end{ex}

\

By making use of Lemma \ref{lem2} with $\varepsilon=0$ and applying Lemma \ref{lem3}, we readily obtain the next theorem.

\begin{thm}
If $f(z)\in \mathcal{H}$ is locally univalent in $\mathbb{U}$ and satisfies\\
\\
${\displaystyle 
\sum\limits_{n=2}^{\infty}\left[\left|\sum\limits_{k=1}^{n}\left\{\sum\limits_{j=1}^{k}(-1)^{k-j}j(j+1)\left(
\begin{array}{c}
\alpha\\
k-j
\end{array}
\right)a_j\right\}\left(
\begin{array}{c}
\beta\\
n-k
\end{array}
\right)\right|\right.
}$
\begin{flushright}
${\displaystyle 
+\left.\left|\sum\limits_{k=1}^{n}\left\{\sum\limits_{j=1}^{k}(-1)^{k-j}j(j-1)\left(
\begin{array}{c}
\alpha\\
k-j
\end{array}
\right)a_j\right\}\left(
\begin{array}{c}
\beta\\
n-k
\end{array}
\right)\right|\right]\leqq 2
}$
\end{flushright}
for some real numbers $\alpha$ and $\beta$, then $f(z)\in \mathcal{C}_{\mathcal{H}}$.
\end{thm}

\

Putting $\alpha=\beta=0$ in the above theorem, we arrive at the following result due to Jahangiri and Silverman \cite{JS}.

\begin{thm}
If $f(z)\in \mathcal{H}$ is locally univalent in $\mathbb{U}$ with
$$
\sum\limits_{n=2}^{\infty}n^2|a_n|\leqq 1,
$$
then $f(z)\in \mathcal{C}_{\mathcal{H}}$.
\end{thm}

\

Furthermore, taking $\alpha=1$ and $\beta=0$ in the theorem, we have

\begin{cor} \label{cor}
If $f(z)\in \mathcal{H}$ is locally univalent in $\mathbb{U}$ and satisfies
$$
\sum\limits_{n=2}^{\infty}\left\{n\left|(n+1)a_n-(n-1)a_{n-1}\right|+(n-1)\left|na_n-(n-2)a_{n-1}\right|\right\}\leqq 2,
$$
then $f(z)\in \mathcal{C}_{\mathcal{H}}$.
\end{cor}

\

\begin{ex}
The function
$$
f(z)=-\int_{0}^{z}\dfrac{\log(1-t)}{t}dt+\overline{\Bigl(z+(1-z)\log(1-z)\Bigr)}=z+\sum\limits_{n=2}^{\infty}\dfrac{1}{n^2}z^n+\sum\limits_{n=2}^{\infty}\dfrac{1}{n(n-1)}\overline{z}^n
$$
satisfies the conditions of Corollary \ref{cor} and belongs to the class $\mathcal{C}_{\mathcal{H}}$.

\begin{figure*}[h]
\begin{center}
\includegraphics[width=6.2cm,clip]{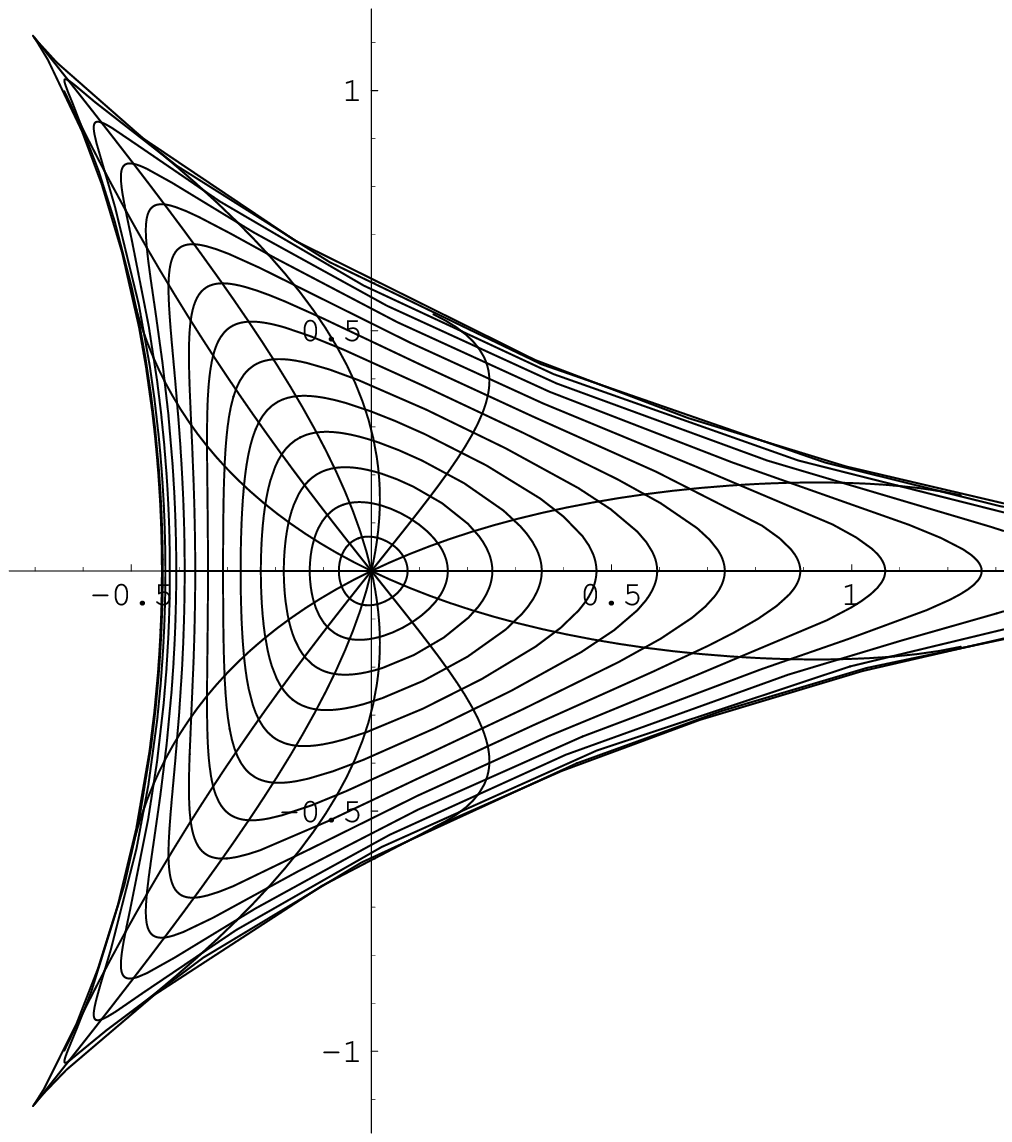}
\end{center}
\caption{The image of ${\displaystyle f(z)=-\int_{0}^{z}\dfrac{\log(1-t)}{t}dt+\overline{\Bigl(z+(1-z)\log(1-z)\Bigr)}}$.}
\end{figure*}
\end{ex}

\

\section{Appendix}

\

A sequence $\left\{c_n\right\}_{n=0}^{\infty}$ of non-negative real numbers is called a convex null sequence if $c_n\rightarrow 0$ as $n\rightarrow \infty$ and
$$
c_n-c_{n+1}\geqq c_{n+1}-c_{n+2}\geqq 0
$$
for all $n$ $(n=0,1,2,\cdots)$.

\

The next lemma was obtained by Fej\'{e}r \cite{F}.

\begin{lem} \label{lem4}
Let $\{c_n\}_{k=0}^{\infty}$ be a convex null sequence. Then, the function
$$
p(z)=\dfrac{c_0}{2}+\sum\limits_{n=1}^{\infty}c_n z^n
$$
is analytic and satisfies ${\rm Re}(p(z))>0$ in $\mathbb{U}$.
\end{lem}

\

Applying the above lemma, we deduce

\begin{thm} \label{thm4}
For some $b$ $(|b|<1)$ and some convex null sequence $\{c_n\}_{n=0}^{\infty}$ with $c_0=2$, the function
$$
f(z)=h(z)+\overline{g(z)}=z+\sum\limits_{n=2}^{\infty}\dfrac{c_{n-1}}{n} z^n+\overline{b\left(z+\sum\limits_{n=2}^{\infty}\dfrac{c_{n-1}}{n} z^n\right)}
$$
belongs to the class $\mathcal{C}_{\mathcal{H}}$.
\end{thm}

\

\begin{proof}
Let us define $F(z)$ by
$$
F(z)=\dfrac{h(z)+\varepsilon g(z)}{1+\varepsilon b}=z+\sum\limits_{n=2}^{\infty}\dfrac{c_{n-1}}{n}z^n
$$
for each $\varepsilon$ $(|\varepsilon|=1)$. Then, we know that
$$
F'(z)=\dfrac{c_0}{2}+\sum\limits_{n=1}^{\infty}c_n z^n\qquad (c_0=2).
$$
By virtue of Lemma \ref{lem1} and Lemma \ref{lem4}, it follows that ${\rm Re}(F'(z))>0$ $(z\in \mathbb{U})$, that is, $F(z)\in \mathcal{C}$. Thus, we conclude that $f(z)=h(z)+\overline{g(z)}\in \mathcal{C}_{\mathcal{H}}$.
\end{proof}

\

In the same manner, we also have

\begin{thm} \label{thm5}
For some $b$ $(|b|<1)$ and some convex null sequence $\{c_n\}_{n=0}^{\infty}$ with $c_0=2$, the function
$$
f(z)=h(z)+\overline{g(z)}=z+\sum\limits_{n=2}^{\infty}\dfrac{1}{n}\left(1+\sum\limits_{j=1}^{n-1}c_j\right)z^n+\overline{b\left(z+\sum\limits_{n=2}^{\infty}\dfrac{1}{n}\left(1+\sum\limits_{j=1}^{n-1}c_j\right)z^n\right)}
$$
belongs to the class $\mathcal{C}_{\mathcal{H}}$.
\end{thm}

\

\begin{proof}
Let us define $F(z)$ by
$$
F(z)=\dfrac{h(z)+\varepsilon g(z)}{1+\varepsilon b}=z+\sum\limits_{n=2}^{\infty}\dfrac{1}{n}\left(1+\sum\limits_{j=1}^{n-1}c_j\right)z^n
$$
for each $\varepsilon$ $(|\varepsilon|=1)$. Then, we know that
$$
(1-z)F'(z)=\dfrac{c_0}{2}+\sum\limits_{n=1}^{\infty}c_n z^n\qquad (c_0=2).
$$
Therefore, by the help of Lemma \ref{lem1} and Lemma \ref{lem4}, we obtain that ${\rm Re}\left((1-z)F'(z)\right)>0$ $(z\in \mathbb{U})$, that is, $F(z)\in \mathcal{C}$ which implies that $f(z)=h(z)+\overline{g(z)}\in \mathcal{C}_{\mathcal{H}}$.
\end{proof}

\

\begin{remark}
The sequence
$$
\{c_n\}_{n=0}^{\infty}=\left\{2,1,\dfrac{2}{3},\cdots,\dfrac{2}{n+1},\cdots\right\}
$$
is a convex null sequence because
$$
\lim_{n\rightarrow \infty}c_n=\lim_{n\rightarrow \infty}\left(\dfrac{2}{n+1}\right)=0,\qquad c_n-c_{n+1}=\dfrac{2}{(n+1)(n+2)}\geqq 0
$$
and
$$
(c_n-c_{n+1})-(c_{n+1}-c_{n+2})=\dfrac{4}{(n+1)(n+2)(n+3)}\geqq 0\qquad (n=0,1,2,\cdots).
$$
\end{remark}

\

Setting $b=\dfrac{1}{4}$ in Theorem \ref{thm4} with the above sequence $\{c_n\}_{n=0}^{\infty}$, we derive

\begin{ex} \label{ex2}
The function
$$
f(z)=-z-2\int_{0}^{z}\dfrac{\log(1-t)}{t}dt-\overline{\dfrac{1}{4}\left(z+2\int_{0}^{z}\dfrac{\log(1-t)}{t}dt\right)}=z+\sum\limits_{n=2}^{\infty}\dfrac{2}{n^2}z^n+\overline{\dfrac{1}{4}\left(z+\sum\limits_{n=2}^{\infty}\dfrac{2}{n^2}z^n\right)}
$$
is in the class $\mathcal{C}_{\mathcal{H}}$.

\begin{figure*}[h]
\begin{center}
\includegraphics[width=10cm,clip]{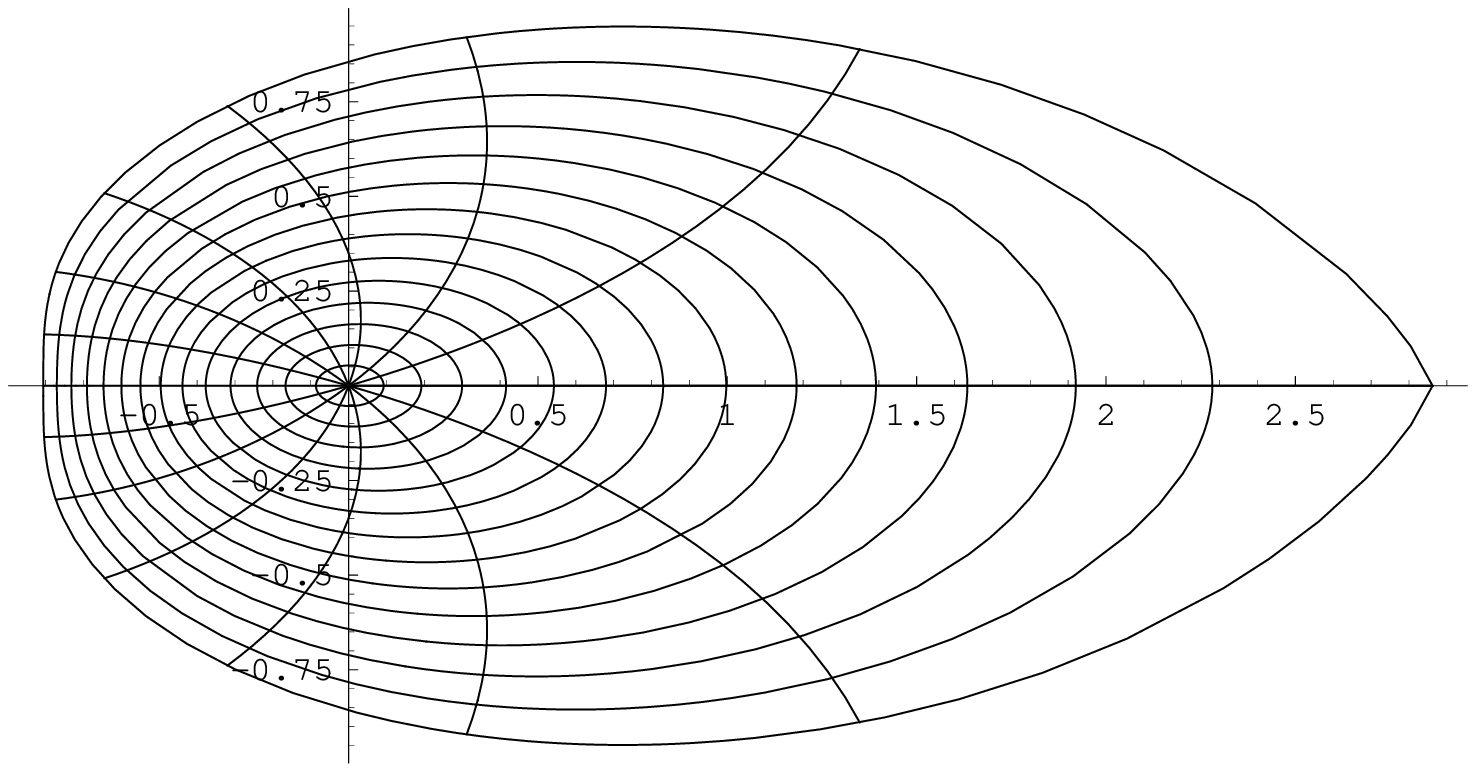}
\end{center}
\caption{The image of $f(z)$ in Example \ref{ex2}.}
\end{figure*}
\end{ex}

\

Moreover, we know

\begin{remark}
The sequence
$$
\{c_n\}_{n=0}^{\infty}=\left\{2,1,\dfrac{1}{2},\cdots,2^{1-n},\cdots\right\}
$$
is a convex null sequence because
$$
\lim_{n\rightarrow \infty}c_n=\lim_{n\rightarrow \infty}2^{1-n}=0,\qquad c_n-c_{n+1}=2^{-n}\geqq 0
$$
and
$$
(c_n-c_{n+1})-(c_{n+1}-c_{n+2})=2^{-(n+1)}\geqq 0\qquad (n=0,1,2,\cdots).
$$
\end{remark}

\

Hence, letting $b=\dfrac{1}{4}$ in Theorem \ref{thm5} with the sequence $\{c_n\}_{n=0}^{\infty}=\{2^{1-n}\}_{n=0}^{\infty}$, we have

\begin{ex} \label{ex3}
The function
\begin{eqnarray*}
f(z) &=& -3\log(1-z)+4\log\left(1-\dfrac{z}{2}\right)+\overline{\left(-\dfrac{3}{4}\log(1-z)+\log\left(1-\dfrac{z}{2}\right)\right)}\\
 & & \\
 &=& z+\sum\limits_{n=2}^{\infty}\dfrac{1}{n}\left(1+\sum\limits_{j=1}^{n-1}2^{1-j}\right)z^n+\overline{\dfrac{1}{4}\left(z+\sum\limits_{n=2}^{\infty}\dfrac{1}{n}\left(1+\sum\limits_{j=1}^{n-1}2^{1-j}\right)z^n\right)}
\end{eqnarray*}
is in the class $\mathcal{C}_{\mathcal{H}}$.

\

\begin{figure*}[h]
\begin{center}
\includegraphics[width=9cm,clip]{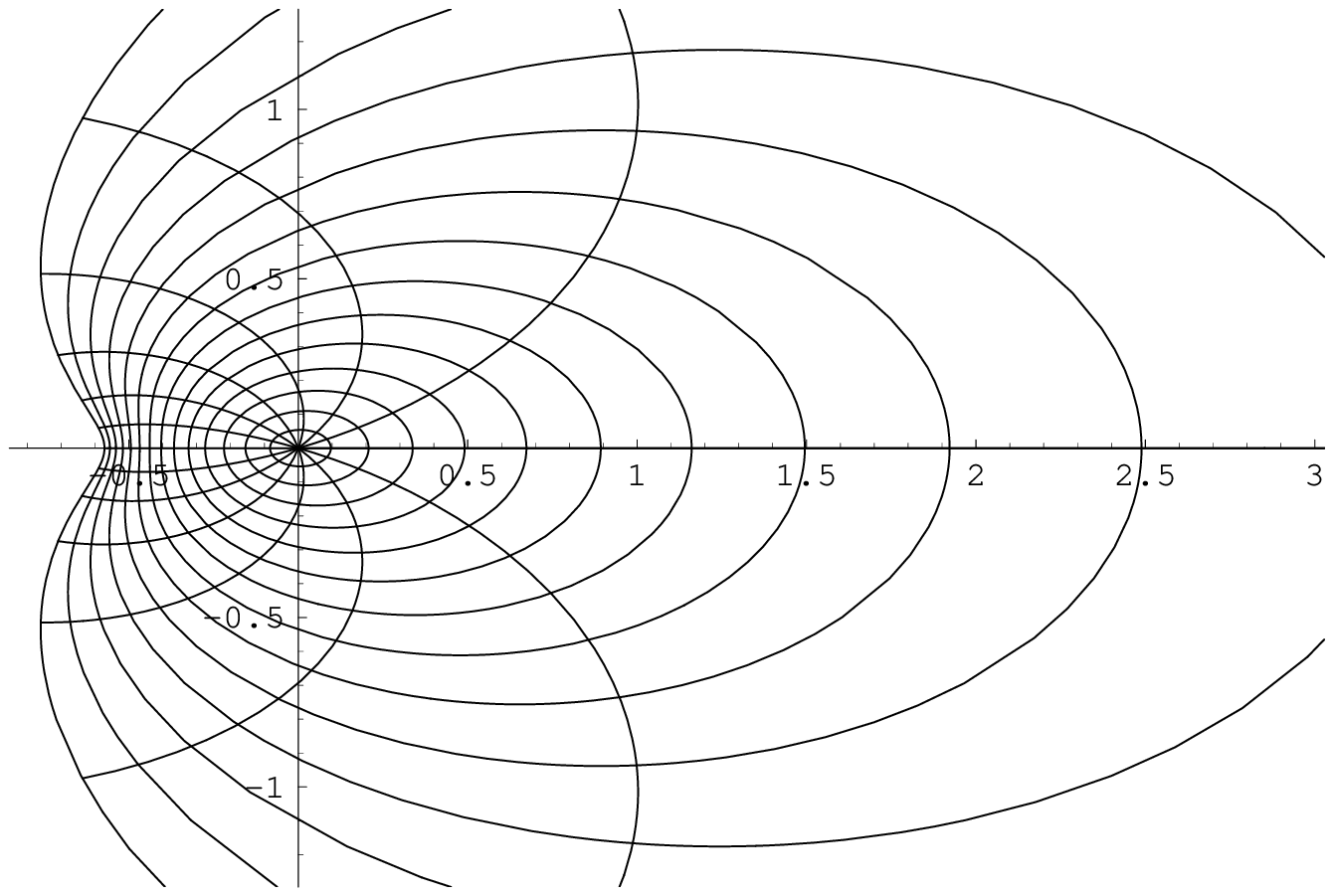}
\end{center}
\caption{The image of $f(z)$ in Example \ref{ex3}.}
\end{figure*}
\end{ex}

\

\end{document}